\documentclass[11pt]{article}
\usepackage{tikz}
\usepackage[all]{xy}
\usepackage{amsfonts}
\usepackage{amsfonts}
\usepackage{amsfonts}
\usepackage{mathrsfs}
\usepackage{bm}
\usepackage{cite}
\usepackage[colorlinks,linkcolor=blue,citecolor=blue]{hyperref}
\usepackage{amssymb,amsmath} 
\usepackage{makecell}
\usepackage{tikz}
 \allowdisplaybreaks
\usetikzlibrary{arrows,shapes,chains}
 \allowdisplaybreaks

\catcode`!=11
\let\!int\int \def\int{\displaystyle\!int}
\let\!lim\lim \def\lim{\displaystyle\!lim}
\let\!sum\sum \def\sum{\displaystyle\!sum}
\let\!sup\sup \def\sup{\displaystyle\!sup}
\let\!inf\inf \def\inf{\displaystyle\!inf}
\let\!cap\cap \def\cap{\displaystyle\!cap}
\let\!max\max \def\max{\displaystyle\!max}
\let\!min\min \def\min{\displaystyle\!min}
\let\!frac\frac \def\frac{\displaystyle\!frac}
\catcode`!=12

\let\oldsection\section
\renewcommand\section{\setcounter{equation}{0}\oldsection}

\allowdisplaybreaks
\def\pf{\it{Proof.}\rm\quad}

\def\N{\mathbb{N}}\def\Z{\mathbb{Z}}
\def\Q{\mathbb{Q}}

\def\t{\widetilde{t}}
\def\S{\widetilde{S}}
\def\z{\zeta}
\def\ZZ{\mathcal{Z}}

\newtheorem{thm}{Theorem}[section]
\newtheorem{lem}[thm]{Lemma}
\newtheorem{cor}[thm]{Corollary}

\newtheorem{exa}{Example}[section]

\begin{document}
\title {\bf Two Variants of Euler Sums}
\author{
{ Ce Xu$^{a,b,}$\thanks{Email: cexu2020@outlook.com}\quad Weiping Wang$^{c,}$\thanks{Email:  wpingwang@yahoo.com, wpingwang@zstu.edu.cn}}\\[1mm]
\small a. School of Mathematics and Statistics, Anhui Normal University,\\ \small Wuhu 241000, P.R. China\\
\small b. Multiple Zeta Research Center, Kyushu University \\
\small  Motooka, Nishi-ku, Fukuoka 819-0389, Japan\\
\small c. School of Science, Zhejiang Sci-Tech University,\\
\small  Hangzhou 310018, P.R. China}

\date{}
\maketitle \noindent{\bf Abstract} For positive integers $p_1,p_2,\ldots,p_k,q$ with $q>1$, we define the Euler $T$-sum $T_{p_1p_2\cdots p_k,q}$ as the sum of those terms of the usual infinite series for the classical Euler sum $S_{p_1p_2\cdots p_k,q}$ with odd denominators. Like the Euler sums, the Euler $T$-sums can be evaluated according to the contour integration and residue theorem. Using this fact, we obtain explicit formulas for Euler $T$-sums with repeated arguments analogous to those known for Euler sums. Euler $T$-sums can be written as rational linear combinations of the Hoffman $t$-values. Using known results for Hoffman $t$-values, we obtain some examples of Euler $T$-sums in terms of (alternating) multiple zeta values. Moreover, we prove an explicit formula of triple $t$-values in terms of zeta values, double zeta values and double $t$-values. We also define alternating Euler $T$-sums and prove some results about them by the contour integration and residue theorem. Furthermore, we define another Euler type $T$-sums and find many interesting results. In particular, we give an explicit formulas of triple Kaneko-Tsumura $T$-values of even weight in terms of single and the double $T$-values. Finally, we prove a duality formula of Kaneko-Tsumura's conjecture.
\\[2mm]


\tableofcontents
\section{Introduction and Notations}

For positive integers $s_1,\ldots,s_k$ with $s_1>1$, the multiple zeta value (MZV for short) is defined by
\begin{align}\label{1.1}
\zeta(s_1,s_2,\ldots,s_k):=\sum\limits_{n_1>\cdots>n_k\geq 1} \frac{1}{n_1^{s_1}n_2^{s_2}\cdots n_k^{s_k}}.
\end{align}
The study of multiple zeta values began in the early 1990s with the works of Hoffman \cite{H1992} and Zagier \cite{DZ1994}. The study of multiple zeta values have attracted a lot of
research in the area in the last two decades. For detailed history and applications, please see the book of Zhao \cite{Z2016}.

In a recent paper \cite{H1}, Hoffman introduced and studied a new kind of multiple zeta values
\begin{align}\label{1.2}
t(s_1,s_2,\ldots,s_k):&=\sum\limits_{n_1>\cdots>n_k\geq 1\atop n_i \ {\rm odd}} \frac{1}{n_1^{s_1}n_2^{s_2}\cdots n_k^{s_k}}\nonumber\\&=\sum\limits_{n_1>\cdots>n_k\geq 1} \frac{1}{(2n_1-1)^{s_1}(2n_2-1)^{s_2}\cdots (2n_k-1)^{s_k}}
\end{align}
which is called multiple $t$-values. As its normalized version,
\begin{align}\label{1.3}
\widetilde{t}(s_1,s_2,\ldots,s_k):&\nonumber=2^{s_1+\cdots+s_k}t(s_1,s_2,\ldots,s_k)\\&=\sum\limits_{n_1>\cdots>n_k\geq 1} \frac{1}{(n_1-1/2)^{s_1}(n_2-1/2)^{s_2}\cdots (n_k-1/2)^{s_k}}
\end{align}
we call it multiple $\widetilde{t}$-values. In both these definitions, we call $k$ the ``depth" and $s_1+\cdots+s_k$ the ``weight".

In this paper we consider the odd variant of Euler sums
\begin{align}\label{1.4}
T_{p_1p_2\cdots p_k,q}:=\sum\limits_{n=1}^\infty \frac{h_{n-1}^{(p_1)}h_{n-1}^{(p_2)}\cdots h_{n-1}^{(p_k)}}{(n-1/2)^q},
\end{align}
which we call Euler $T$-sums, where $p_j\in \N\ (j=1,2,\ldots,k)$ and $2 \leq q\in \N$ with $p_1\leq p_2\leq \cdots \leq p_k$. Here  $h_n^{(p)}$ is defined for $n\in\N_0,p\in\N$ by
\begin{align}
h_n^{(p)}:=\sum\limits_{k=1}^n \frac{1}{(k-1/2)^p},\quad h_0^{(p)}:=0,\quad h_n:=h_n^{(1)}.
\end{align}
The classical Euler sum was introduced by Flajolet and Salvy \cite{FS1998}, which is defined by
\begin{align}\label{1.6}
S_{p_1p_2\cdots p_k,q}:=\sum\limits_{n=1}^\infty \frac{H_{n}^{(p_1)}H_{n}^{(p_2)}\cdots H_{n}^{(p_k)}}{n^q},
\end{align}
where $H_n^{(p)}$ is harmonic number of order $p$ defined by
\begin{align}
H_n^{(p)}:=\sum\limits_{k=1}^n \frac{1}{n^p},\quad H_0^{(p)}:=0,\quad H_n:=H_n^{(1)}.
\end{align}
In the definitions of (\ref{1.4}) and (\ref{1.6}), the quantity $s_1+\cdots+s_k+q$ is called the ``weight" of the sum, and the quantity $k$ is called the ``degree". The linear sums $S_{p,q}$ was the first considered by Euler in 1742 (see \cite{B1989} for a discussion). Classical Euler sums may be studied through a profusion of methods: combinatorial, analytic and algebraic. There are many other researches on Euler sums and Euler type sums. Some related results for Euler sums may be seen in the works of \cite{BBG1994,BBG1995,BZB2008,F2005,M2014,W2017} and references therein.

Since repeated summands in partitions are indicated by powers, we denote, for instance, the sum
\begin{align*}
T_{1^22^35,q}:=\sum\limits_{n=1}^\infty \frac{h_{n-1}^{2}\left(h_{n-1}^{(2)}\right)^3 h_{n-1}^{(5)}}{(n-1/2)^q}.
\end{align*}
As remarked by Flajolet and Salvy \cite{FS1998}, every Euler sum of weight $w$ and degree $k$ is a $\mathbb{Q}$-linear
combination of MZVs of weight $w$ and depth at most $k + 1$ (explicit formula see our previous paper \cite{XW2018}). According to the definitions of Hoffman $t$-value and Euler $T$-sum, it is clear that every Euler $T$-sum of weight $w$ and degree $k$ is a $\mathbb{Q}$-linear
combination of Hoffman $t$-value of weight $w$ and depth at most $k + 1$. Because, by the methods of \cite{XW2018,H2016}, we may easily deduce the following relation
\begin{align*}
T_{i_1i_2\cdots i_m,q}
    =\sum_{\xi\in\mathcal{C}_m}\sum_{\sigma\in\mathcal{S}_m}
        \frac{\t(q,J_1(I_{\sigma}^{(m)}),J_2(I_{\sigma}^{(m)}),\ldots,J_p(I_{\sigma}^{(m)}))}{\xi_1!\xi_2!\cdots\xi_p!}
       ,
\end{align*}
where $\xi:=(\xi_1,\xi_2,\ldots,\xi_p)\in\mathcal{C}_m$ ($\mathcal{C}_m$ is a set of all compositions of $m$) and a permutation $\sigma\in\mathcal{S}_m$ ($\mathcal{S}_m$ is a symmetric group of all permutations on $m$ symbols), $I_{\sigma}^{(m)}=(i_{\sigma(1)},\ldots,i_{\sigma(m)})$, and
\[
J_c(I_{\sigma}^{(m)})=i_{\sigma(\xi_1+\cdots+\xi_{c-1}+1)}+\cdots+i_{\sigma(\xi_1+\cdots+\xi_c)}\,,
    \quad\text{for } c=1,2,\ldots,p\,.
\]

The motivation for this paper arises from the results of Flajolet and Salvy. In \cite{FS1998}, Flajolet and Salvy used the method of contour integration to evaluated the classical Euler sums $S_{p_1p_2\cdots p_k,q}$.
Contour integration is a classical technique for evaluating infinite sums by reducing them to a finite number of residue computations. They used the method to found many interesting results. In particular, they proved the famous result that a nonlinear Euler sum $S_{p_1p_2\cdots p_k,q}$ reduces to a combination of sums of lower orders whenever the weight $p_1+p_2+\cdots+p_k+q$ and the order $k$ are of the same parity.
In this paper, we will extend their method to Euler $T$-sums and find many similar results.

The main purpose of this paper is study the reducible formulas of Euler $T$-sums and type $T$-sums by the method of Contour integral. We will prove that a nonlinear Euler $T$-sum $T_{p_1p_2\cdots p_k,q}$ reduces to a combination of $\log(2)$, Euler $T$-sums with depth $\leq k-1$, multiple zeta values with depth $\leq k$ whenever the weight $p_1+p_2+\cdots+p_k+q$ and the order $k$ are of the same parity.

The remainder of this paper is organized as follows. In the second section we provide some asymptotic formulas of $\Psi(1/2-s)$. Then we apply it and contour integral to evaluate the linear and nonlinear Euler $T$-sums. Specially, we establish the explicit formulas of linear $T$-sum $T_{p,q}$ with $p+q$ odd, quadratic $T$-sums $T_{p_1p_2,q}$ with $p_1+p_2+q$ even and cubic $T$-sum $T_{1^3,q}$ with $q$ even. Further, we prove that all Euler $T$-sum $T_{p_1p_2\cdots p_k,q}$ can be expressed in terms of a combination of $\log(2)$, Euler $T$-sums with degree $\leq k-1$, multiple zeta values with depth $\leq k$ whenever the weight $p_1+p_2+\cdots+p_k+q$ if the weight $p_1+p_2+\cdots+p_k+q$ and order $k$ are of the same parity. In the third section, we define an alternating Euler $T$-sums and evaluate the linear and a quadratic alternating Euler $T$-sums. In the fourth section, we define an Euler type $T$-sums $\S_{p_1p_2\cdots p_k,q}$, which is defined by
\begin{align}
\S_{p_1p_2\cdots p_k,q}:=\sum\limits_{n=1}^\infty \frac{h_n^{(p_1)}h_n^{(p_2)}\cdots h_n^{(p_k)}}{n^q}.
\end{align}
Then, we establish many relations of the sum by using the method of contour integration. In particular, we prove a general formula of quadratic sums $\S_{p_1p_2,q}$ with $p_1+p_2+q$ even. According to the relation of $\S_{p_1p_2,q}$ and triple Kaneko-Tsumura $T$-value, we can obtain a formula of triple $T$-value with weight even. In the last section, we prove a duality identity of Kaneko-Tsumura's conjecture, and establish a relation between the double $T$-values and the double $t$-values.

\section{Evaluations of Euler $T$-sums}

In \cite{XC2019}, the first author defined a parametric digamma (or Psi) function $\Psi (-s;a)$ by
\begin{align}\label{2.1}
\Psi \left( { - s;a} \right) + \gamma : = \frac{1}{{s - a}} + \sum\limits_{k = 1}^\infty  {\left( {\frac{1}{{k + a}} - \frac{1}{{k + a - s}}} \right)} ,\quad (s\in \mathbb{C},\ a\in \mathbb{C}\setminus\N^-).
\end{align}
The function $\Psi (-s;a)$ is meromorphic in the entire complex plane with a simple pole at $s=n+a$ for each negative integer $n$. In here, we let
$$\Psi(-s):=\Psi(-s;-1/2)+\gamma=\frac{1}{s+1/2}+\sum\limits_{k=1}^\infty \left(\frac{1}{k-1/2}-\frac{1}{k-1/2-s} \right).$$
From Theorems 1.1-1.3 and Corollary 2.4 in \cite{XC2019}, by direct calculations we can obtain the following identities $(2\leq p\in\N)$
\begin{align}
&\begin{aligned}\Psi \left( { \frac{1}{2}- s} \right)\mathop  = \limits^{s \to  n}& \frac{1}{s-n}+H_n+2\log(2)\\&+\sum\limits_{j=1}^\infty \left((-1)^j H_n^{(j+1)}-\zeta(j+1) \right)(s-n)^j\quad (n\in\N_0:=\N\cup \{0\}),\end{aligned}\label{2.2}\\
&\begin{aligned}\Psi \left( { \frac{1}{2}- s} \right)\mathop  = \limits^{s \to  n-1/2}&h_n+\sum\limits_{j=1}^\infty \left((-1)^jh_n^{(j+1)}-\widetilde{t}(j+1)\right)(s+1/2-n)^j\quad(n\in\N_0),\end{aligned}\label{2.3}\\
&\begin{aligned}\Psi \left( { \frac{1}{2}- s} \right)\mathop  = \limits^{s \to  -(n-1/2)}&h_{n-1}+\sum\limits_{j=1}^\infty \left(h_{n-1}^{(j+1)}-\widetilde{t}(j+1)\right)(s-1/2+n)^j\quad (n \in \N),\end{aligned}\label{2.4}\\
&\begin{aligned}\frac{{{\Psi ^{\left( {p - 1} \right)}}\left( { \frac1{2}- s} \right)}}{{\left( {p - 1} \right)!}}\mathop  = \limits^{s \to n } \frac{1}{{{{\left( {s - n } \right)}^p}}}\left(1
+ (-1)^p\sum\limits_{j = p}^\infty  \binom{j-1}{p-1} \left( {\zeta \left( {j} \right) + {{\left( { - 1} \right)}^j}H_n^{\left( j \right)}} \right){\left( {s - n } \right)^j}\right)\quad (n\in\N_0),\end{aligned}\label{2.5}\\
&\begin{aligned}\frac{{{\Psi ^{\left( {p - 1} \right)}}\left( { \frac1{2}- s} \right)}}{{\left( {p - 1} \right)!}}\mathop  = \limits^{s \to n-1/2 } (-1)^p
\sum\limits_{j = p}^\infty  \binom{j-1}{p-1} \left( \widetilde{t}(j)+(-1)^jh_n^{(j)} \right){\left( {s - n+1/2 } \right)^{j-p}}\quad (n\in\N_0),\label{2.6}\end{aligned}\\
&\begin{aligned}\frac{{{\Psi ^{\left( {p - 1} \right)}}\left( { \frac1{2}- s} \right)}}{{\left( {p - 1} \right)!}}\mathop  = \limits^{s \to -(n-1/2) } (-1)^p
\sum\limits_{j = p}^\infty  \binom{j-1}{p-1} \left( \widetilde{t}(j)-h_{n-1}^{(j)} \right){\left( {s + n-1/2 } \right)^{j-p}}\quad (n\in\N).\label{2.7}\end{aligned}
\end{align}
We also deduce that from \cite{XC2019}
\begin{align}
&\pi \tan \left( {\pi s} \right)\mathop  = \limits^{s \to (n-1/2)}  - \frac{1}
{{s - \frac{{2n - 1}}
{2}}} + 2\sum\limits_{k = 1}^\infty  {\zeta \left( {2k} \right){{\left( {s - \frac{{2n - 1}}
{2}} \right)}^{2k - 1}}} \quad (n\in \Z).\label{2.8}
\end{align}

We define a kernel function $\xi \left( s \right)$ by the two requirements: 1. $\xi \left( s \right)$ is meromorphic in the whole complex plane. 2. $\xi \left( s \right)$ satisfies $\xi \left( s \right)=o(s)$ over an infinite collection of circles $\left| s \right| = {\rho _k}$ with ${\rho _k} \to \infty . $

\begin{lem}(\cite{FS1998})\label{lem1}
Let $\xi \left( s \right)$ be a kernel function and let $r(s)$ be a rational function which is $O(s^{-2})$ at infinity. Then
\begin{align}\label{3.1}
\sum\limits_{\alpha  \in O} {{\mathop{\rm Res}}{{\left[ {r\left( s \right)\xi \left( s \right)},s = \alpha  \right]}}}  + \sum\limits_{\beta  \in S}  {{\mathop{\rm Res}}{{\left[ {r\left( s \right)\xi \left( s \right)},s = \beta  \right]}}}  = 0.
\end{align}
where $S$ is the set of poles of $r(s)$ and $O$ is the set of poles of $\xi \left( s \right)$ that are not poles $r(s)$ . Here ${\rm Res}[r(s),s=\alpha] $ denotes the residue of $r(s)$ at $s= \alpha$.
\end{lem}

In below, we use the identities (\ref{2.2})-(\ref{2.8}) and residue theorem to evaluate some Euler $T$-sums. We need the formula (\cite{A2000,G1989})
\begin{align}
&\pi \tan(\pi s)=2\sum\limits_{k=1}^\infty \widetilde{t}(2k)s^{2k-1}=\sum\limits_{k=1}^\infty \frac{(-1)^{k-1}2^{2k}(2^{2k}-1)B_{2k}\pi^{2k}}{(2k)!}s^{2k-1},\label{2.10}
\end{align}
where $B_{2k}$ is Bernoulli numbers. By (\ref{2.10}), we have
\begin{align}
\lim\limits_{s\rightarrow n} \frac{d^p}{ds^p}(\pi\tan(\pi s))&=\lim\limits_{t\rightarrow 0} \frac{d^p}{dt^p}(\pi\tan(\pi t))=\lim\limits_{t\rightarrow 0} \frac{d^p}{dt^p}\left(2\sum\limits_{k=1}^\infty \widetilde{t}(2k)t^{2k-1}\right)\nonumber\\
&=(1-(-1)^p)p!\t(p+1).\label{2.11}
\end{align}
Next, let ${\rm{Res}}\left[ {f\left( s \right),s = \alpha } \right]$ to denote the residue of $f(s)$ at $s=\alpha$.
\begin{lem}\label{lem2} If a meromorphic function $F(s)$ has pole of order $m$ at $s=\alpha$, then
\begin{align*}
{\rm Res}[F(s),s=\alpha]=\lim\limits_{s\rightarrow\alpha}\frac{1}{(m-1)!}\frac{d^{m-1}}{ds^{m-1}}\left[(s-\alpha)^m F(s)\right]=\lim\limits_{s\rightarrow\alpha}\frac{1}{m!}\frac{d^{m}}{ds^{m}}\left[(s-\alpha)^{m+1} F(s)\right].
\end{align*}
\end{lem}
\pf This proof is very simple, so we omitted.\hfill$\square$

\subsection{Linear Euler $T$-sums}
In this subsection, we will prove the linear $T$-sum $T_{p,q}$ can be expressed in terms of $\log(2)$, zeta values and $\widetilde{t}$-values with $p+q$ odd $(q\geq 2)$.

\begin{thm}\label{thm1} For positive integer $q>1$,
\begin{align}\label{2.12}
(1+(-1)^q)T_{1,q}=&(-1)^{q+1}\widetilde{t}(q+1)+(1+(-1)^q)2\log(2)\widetilde{t}(q)\nonumber\\&-2\sum\limits_{2k_1+k_2=q,\atop k_1,k_2\geq 1} \widetilde{t}(2k_1)\zeta(k_2+1).
\end{align}
\end{thm}
\begin{thm}\label{thm2} For positive integers $p,q>1$,
\begin{align}\label{2.13}
(1-(-1)^{p+q})T_{p,q}=&(-1)^{p+q}\widetilde{t}(p+q)-(-1)^p(1+(-1)^q)\widetilde{t}(p)\widetilde{t}(q)\nonumber\\
&-(-1)^{p}\sum\limits_{k=0}^{p-1}((-1)^k-1)\binom{p+q-k-2}{q-1} \widetilde{t}(k+1) \zeta(p+q-k-1)\nonumber\\
&+2(-1)^p\sum\limits_{2k_1+k_2=q+1,\atop k_1,k_2\geq 1} \binom{k_2+p-2}{p-1} \t(2k_1)\z(k_2+p-1).
\end{align}
\end{thm}
\pf In the context of this paper, these theorems results form applying the kernels
\begin{align*}
{\pi\tan(\pi s)\Psi(1/2-s)}\quad{\rm and}\quad \frac{\pi \tan(\pi s)\Psi^{(p-1)}(1/2-s)}{(p-1)!}
\end{align*}
to the base function $r(s)=s^{-q}$, respectively. Now, we only prove the formula (\ref{2.13}). The identity (\ref{2.12}) can be shown in a similar way. Let
\begin{align*}
F_{p-1,q}(s):=\frac{\pi \tan(\pi s)\Psi^{(p-1)}(1/2-s)}{s^q(p-1)!}.
\end{align*}
The function $F_{p-1,q}(s)$ only have poles at the $s=0,\pm (n-1/2),n\ (n\in\N)$. At a positive integer $n$, the pole $\pm(n-1/2)$ are simple and the residue is
\begin{align*}
&{\rm Res}[F_{p-1,q}(s),s=n-1/2]=-\lim\limits_{s\rightarrow n-1/2} \frac{\Psi^{(p-1)}(1/2-s)}{s^q(p-1)!}=-\frac{(-1)^p\t(p)+h_n^{(p)}}{(n-1/2)^q},\\
&{\rm Res}[F_{p-1,q}(s),s=1/2-n]=-\lim\limits_{s\rightarrow 1/2-n} \frac{\Psi^{(p-1)}(1/2-s)}{s^q(p-1)!}=-(-1)^{p+q}\frac{\t(p)-h_{n-1}^{(p)}}{(n-1/2)^q},
\end{align*}
where we used the identities (\ref{2.6}) and (\ref{2.7}). For a positive intger $n$, the pole has order $p-1$ (since $s=n$ is a zero of order one of $\tan(\pi s)$), by (\ref{2.5}), (\ref{2.11}) and Lemma \ref{lem2}, the residue is
\begin{align*}
{\rm Res}[F_{p-1,q}(s),s=n]&=\frac{1}{(p-1)!}\lim\limits_{s\rightarrow n}\frac{d^{p-1}}{ds^{p-1}} \frac{\pi \tan(\pi s)}{s^q}\\&=(-1)^{p-1} \sum\limits_{k=0}^{p-1}((-1)^k-1)\binom{p+q-k-2}{q-1} \frac{\t(k+1)}{n^{p+q-k-1}}.
\end{align*}
By (\ref{2.5}) and (\ref{2.10}) with $n=0$, we know that if $s\rightarrow 0$, then
\begin{align*}
F_{p-1,q}(s)=\frac{2}{s^{p+q-1}} \left\{\sum\limits_{k=1}^\infty \t(2k)s^{2k-2}+(-1)^p \sum\limits_{k_1,k_2=1}^\infty \binom{k_2+p-2}{p-1} \t(2k_1)\z(k_2+p-1)s^{2k_1+k_2+p-3} \right\}.
\end{align*}
Hence, the residue of the pole of order $p+q-1$ at $0$ is found to be
\begin{align*}
{\rm Res}[F_{p-1,q}(s),s=0]=(1+(-1)^{p+q})\t(p+q)+2(-1)^p \sum\limits_{2k_1+k_2=q+1,\atop k_1,k_2\geq 1} \binom{k_2+p-2}{p-1} \t(2k_1)\z(k_2+p-1).
\end{align*}
Summing these four contributions yields the statement of the theorem {\ref{thm2}}.\hfill$\square$

Therefore, from Theorems \ref{thm1} and \ref{thm2}, we arrive at the conclusions
\begin{align*}
&\t(q,1)=T_{1,q}\in \Q[\log(2), \text{Zeta values}]\quad (q\ \text{is even}),\\
&\t(q,p)=T_{p,q}\in \Q[\text {Zeta values}]\quad (p+q \text{ is odd}),
\end{align*}
where we used the relaton $\t(p)=(2^p-1)\z(p)$. For even weights, a modified form of the identity holds, but without any linear Euler $T$-sum occurring. This gives back well-known nonlinear relations between $\t$-values at even arguments.

\begin{exa} We have 
\begin{align*}
&T_{1,2}=-\frac{7}{2}\zeta(3)+\pi^2\log(2),\\
&T_{2,3}=-\frac{31}{2}\zeta(5)+\frac3{2}\pi^2\z(3),\\
&T_{3,2}=-\frac{31}{2}\zeta(5)+2\pi^2\z(3),\\
&T_{1,4}=-\frac{31}{2}\zeta(5)+\frac1{3}\pi^4\log(2)-\frac1{2}\pi^2\z(3).
\end{align*}
\end{exa}

\subsection{Quadratic Euler $T$-sums}

\begin{thm}\label{thm2.5} For positive integer $q>1$,
\begin{align}\label{2.14}
(1+(-1)^q)T_{1^2,q}&=\pi^2\zeta(q)+(-1)^q\t(q+2)+(1-(-1)^q)4\log(2)\t(q+1)-2\t(q+1,1)\nonumber\\&\quad+(1+(-1)^q)4\log^2(2)\t(q)-4\sum\limits_{2k_1+k_2=q+1,\atop k_1,k_2\geq 1} \t(2k_1)\z(k_2+1)\nonumber\\
&\quad-8\log(2)\sum\limits_{2k_1+k_2=q,\atop k_1,k_2\geq 1} \t(2k_1)\z(k_2+1)\nonumber\\&\quad+2\sum\limits_{2k_1+k_2+k_3=q,\atop k_1,k_2\geq1} \t(2k_1)\z(k_2+1)\z(k_3+1).
\end{align}
\end{thm}

\pf The proof is based on the function
\begin{align*}
F_{0^2,q}(s):=\frac{\pi \tan(\pi s)\left(\Psi(1/2-s) \right)^2}{s^q}
\end{align*}
and the usual residue computation. The function $F_{0^2,q}(s)$ only have poles at the $s=0,\pm (n-1/2),n\ (n\in\N)$.  By a similar argument as in the proof of (\ref{2.13}), we deduce
\begin{align*}
&{\rm Res}[F_{0^2,q}(s),s=n]=\frac{\pi^2}{n^q}\quad(n\in\N),\\
&{\rm Res}[F_{0^2,q}(s),s=n-1/2]=-\frac{h_n^2}{(n-1/2)^q}\quad(n\in\N),\\
&{\rm Res}[F_{0^2,q}(s),s=1/2-n]=-(-1)^q\frac{h_{n-1}^2}{(n-1/2)^q}\quad(n\in\N)
\end{align*}
and
\begin{align*}
{\rm Res}[F_{0^2,q}(s),s=0]&=(1+(-1)^q)\t(q+2)+(1-(-1)^q)4\log(2)\t(q+1)+(1+(-1)^q)4\log^2(2)\t(q)\\
&\quad-4\sum\limits_{2k_1+k_2=q+1,\atop k_1,k_2\geq 1} \t(2k_1)\z(k_2+1)-8\log(2)\sum\limits_{2k_1+k_2=q,\atop k_1,k_2\geq 1} \t(2k_1)\z(k_2+1)\\
&\quad-2\sum\limits_{2k_1+1k_2+k_3=q,\atop k_1,k_2\geq1} \t(2k_1)\z(k_2+1)\z(k_3+1).
\end{align*}
Thus, summing these four contributions yields the desired result.\hfill$\square$

Hence, from (\ref{2.14}),
\begin{align*}
T_{1^2,q}+\t(q+1,1)\in \Q[\log(2),\text{Zeta values}]\quad (q\ \text{even}).
\end{align*}
If $q=2$, then
$$T_{1^2,2}+\t(3,1)=2\log^2(2)\pi^2.$$
Note that $T_{1^2,q}=2\t(q,1,1)+\t(q,2).$

\begin{thm}\label{thm2.6} For positive integer $q>1$,
\begin{align}\label{2.15}
(1-(-1)^q)T_{12,q}&=\pi^2(\z(q,1)+2\log(2)\z(q))-(q-1)\pi^2\z(q+1)+((-1)^{q+1}-1)\t(2)\t(q,1)\nonumber\\
&\quad-\t(2)\t(q+1)-\t(q+2,1)-\t(q+1,2)-(-1)^q\t(q+3)\nonumber\\
&\quad+2(1+(-1)^q)\log(2)\t(q+2)+2\sum\limits_{2k_1+k_2=q+2,\atop k_1,k_2\geq 1} (k_2-1)\t(2k_1)\z(k_2+1)\nonumber\\
&\quad+4\log(2)\sum\limits_{2k_1+k_2=q+1,\atop k_1,k_2\geq 1} k_2 \t(2k_1)\z(k_2+1)\nonumber\\&\quad-2\sum\limits_{2k_1+k_2+k_3=q+1,\atop k_1,k_2\geq 1} k_3 \t(2k_1)\z(k_2+1)\z(k_3+1)
\end{align}
\end{thm}
\pf Consider the function
\begin{align*}
F_{01,q}(s):=\frac{\pi \tan(\pi s)\Psi(1/2-s)\Psi^{(1)}(1/2-s) }{s^q}.
\end{align*}
Then, by a similar argument as in the proof of (\ref{2.14}), we can prove the theorem.\hfill$\square$

If $q=3$, then
\[2T_{12,3}+\t(5,1)+\t(4,2)=-\frac{\pi^6}{24}+6\pi^2\log(2)\z(3).\]

\begin{thm} \label{thm2.7} For positive integer $p,q>1$,
\begin{align}
&(1-(-1)^{p+q})T_{1p,q}\nonumber\\&=-(-1)^p(1+(-1)^q)\t(p)T_{1,q}-(-1)^p\t(p)\t(q+1)\nonumber\\&\quad-T_{1,p+q}-T_{p,q+1}-(-1)^{p+q}\t(p+q+1)\nonumber\\
&\quad+(1+(-1)^{p+q})2\log(2)\t(p+q)-2\sum\limits_{2k_1+k_2=p+q,\atop k_1,k_2\geq 1} \t(2k_1)\z(k_2+1)\nonumber\\
&\quad+(-1)^p 2 \sum\limits_{2k_1+k_2=q+2,\atop k_1,k_2\geq 1} \binom{k_2+p-2}{p-1} \t(2k_1)\z(k_2+p-1)\nonumber\\
&\quad+(-1)^p4\log(2)\sum\limits_{2k_1+k_2=q+1,\atop k_1,k_2\geq 1} \binom{k_2+p-2}{p-1} \t(2k_1)\z(k_2+p-1)\nonumber\\
&\quad-(-1)^p 2 \sum\limits_{2k_1+k_2+k_3=q+1,\atop k_1,k_2\geq 1}\binom{k_3+p-2}{p-1} \t(2k_1)\z(k_2+1)\z(k_3+p-1)\nonumber\\
&\quad+(-1)^p \sum\limits_{l=0}^p ((-1)^l-1) \binom{p+q-l-1}{q-1} \t(l+1)\z(p+q-l)\nonumber\\
&\quad-(-1)^p  \sum\limits_{l=0}^{p-1} ((-1)^l-1) \binom{p+q-l-2}{q-1}\t(l+1)\nonumber\\
&\quad\quad\quad\quad\times\left( S_{1,p+q-l-1}+2\log(2)\z(p+q-l-1)\right)\nonumber\\
&\quad+(-1)^p \sum\limits_{k=1}^{p-1} (-1)^{k+1} \sum\limits_{l=0}^{p-k-1} ((-1)^l-1) \binom{p+q-k-l-2}{q-1} \t(l+1)\nonumber\\
&\quad\quad\quad\quad\times \left((-1)^kS_{k+1,p+q-k-l-1}-\z(k+1)\z(p+q-k-l-1) \right).
\end{align}
\end{thm}
\pf By computing the residues of the function
\begin{align*}
F_{0(p-1),q}(s):=\frac{\pi \tan(\pi s)\Psi(1/2-s)\Psi^{(p-1)}(1/2-s) }{s^q},
\end{align*}
we may deduce the desired formula. \hfill$\square$

Putting $p=3,q=2$ in equation above gives
\begin{align*}
 2T_{13,2}+T_{1,5}+T_{3,3}=8\log(2)\pi^2\z(3)-\frac 7{360}\pi^6.
\end{align*}
A more general reduction results from the kernel
\[F_{(p_1-1)(p_2-1),q}(s):=\frac{\pi \tan(\pi s)\Psi^{(p_1-1)}(1/2-s)\Psi^{(p_2-1)}(1/2-s)}{(p_1-1)!(p_2-1)!}\]
but it involves a parity restriction on the weight because of its trigonometric factor.

\begin{thm}\label{thm2.8} For positive integer $p_1,p_2,q>1$,
\begin{align*}
(1+(-1)^{p_1+p_2+q})T_{p_1p_2,q}\in \Q[\text{\rm zeta values, double zeta values, double $\t$-values}].
\end{align*}
We have
\begin{align*}
&(1+(-1)^{p_1+p_2+q})T_{p_1p_2,q}\nonumber\\
&=-T_{p_1,p_2+q}-T_{p_2,p_1+q}-(-1)^{p_1+p_2}(1+(-1)^q)\t(p_1)\t(p_2)\t(q)+(-1)^{p_1}((-1)^{p_2+q}-1)\t(p_1)T_{p_2,q}\nonumber\\
&\quad+(-1)^{p_2}((-1)^{p_1+q}-1)\t(p_2)T_{p_1,q}-(-1)^{p_1}\t(p_1)\t(p_2+q)-(-1)^{p_2}\t(p_2)\t(p_1+q)\nonumber\\
&\quad-(-1)^{p_1+p_2}\sum\limits_{l=0}^{p_1+p_2-1}((-1)^l-1)\binom{p_1+p_2+q-l-2}{q-1}\t(l+1)\z(p_1+p_2+q-l-1)\nonumber\\
&\quad+(-1)^{p_1+p_2}\sum\limits_{k=1}^{p_2} (-1)^k \binom{k+p_1-2}{p_1-1} \sum\limits_{l=0}^{p_2-k}((-1)^l-1)\binom{p_2+q-k-l-1}{q-1}\t(l+1)\nonumber\\
&\quad\quad\quad\quad\quad\quad\quad\quad\times \left( \z(k+p_1-1)\z(p_2+q-k-l)+(-1)^{k+p_1-1}S_{k+p_1-1,p_2+q-k-l}\right)
\nonumber\\
&\quad+(-1)^{p_1+p_2}\sum\limits_{k=1}^{p_1} (-1)^k \binom{k+p_2-2}{p_2-1} \sum\limits_{l=0}^{p_1-k}((-1)^l-1)\binom{p_1+q-k-l-1}{q-1}\t(l+1)\nonumber\\
&\quad\quad\quad\quad\quad\quad\quad\quad\times \left( \z(k+p_2-1)\z(p_1+q-k-l)+(-1)^{k+p_2-1}S_{k+p_2-1,p_1+q-k-l}\right)\nonumber\\
&\quad+(-1)^{p_1+p_2+q} \t(p_1+p_2+q) +(-1)^{p_1}2\sum\limits_{k_1+2k_2=p_2+q+1,\atop k_1,k_2\geq 1} \binom{p_1+k_1-2}{p_1-1}\z(p_1+k_1-1)\t(2k_2)\nonumber\\
&\quad+(-1)^{p_2}2\sum\limits_{k_1+2k_2=p_1+q+1,\atop k_1,k_2\geq 1} \binom{p_2+k_1-2}{p_2-1}\z(p_2+k_1-1)\t(2k_2)\nonumber\\
&\quad+(-1)^{p_1+p_2} 2\sum\limits_{k_1+k_2+2k_3=q+2,\atop k_1,k_2,k_3\geq 1} \binom{p_1+k_1-2}{p_1-1}\binom{p_2+k_2-2}{p_2-1}\z(p_1+k_1-1)\z(p_2+k_2-1)\t(2k_3).
\end{align*}
\end{thm}
\pf
Let \[F_{(p_1-1)(p_2-1),q}(s):=\frac{\pi \tan(\pi s)\Psi^{(p_1-1)}(1/2-s)\Psi^{(p_2-1)}(1/2-s)}{s^q(p_1-1)!(p_2-1)!}.\]
By (\ref{2.5})-(\ref{2.8}) and (\ref{2.10}), we arrive at
\begin{align*}
{\rm Res}[F_{(p_1-1)(p_2-1),q}(s),s=n-1/2]&=-\frac{(-1)^{p_1+p_2}\t(p_1)\t(p_2)+(-1)^{p_1}\t(p_1)h_n^{(p_2)}+(-1)^{p_2}\t(p_2)h_n^{(p_1)}}{(n-1/2)^q}\\
&\quad-\frac{h_n^{(p_1)}h_n^{(p_2)}}{(n-1/2)^q},\\
{\rm Res}[F_{(p_1-1)(p_2-1),q}(s),s=1/2-n]&=-(-1)^{p_1+p_2+q}\frac{\t(p_1)\t(p_2)-\t(p_1)h_{n-1}^{(p_2)}-\t(p_2)h_{n-1}^{(p_1)}}{(n-1/2)^q}\\
&\quad-(-1)^{p_1+p_2+q}\frac{h_{n-1}^{(p_1)}h_{n-1}^{(p_2)}}{(n-1/2)^q}
\end{align*}
and
\begin{align*}
&\sum\limits_{n=1}^\infty {\rm Res}[F_{(p_1-1)(p_2-1),q}(s),s=n]\in\Q[\text{zeta values, double zeta values}],\\
&{\rm Res}[F_{(p_1-1)(p_2-1),q}(s),s=0]\in\Q[\text{zeta values}].
\end{align*}
Hence, using Lemma \ref{lem1} and combining the four identities gives the desired result. \hfill$\square$

If $p_1=p_2=q=2$, then
\[T_{2^2,2}+T_{2,4}=\frac{7}{360}\pi^6.\]

By the definitions of $\t$-values and $T$-sums, we have $(p_1\geq 1,p_2,q>1)$
\begin{align*}
T_{p_1p_2,q}&=\t(q,p_1,p_2)+\t(q,p_2,p_1)+\t(q,p_1+p_2)\\
&=\t(p_2)\t(q,p)-\sum\limits_{n=1}^\infty \frac{h_{n-1}^{(p_1)}(\t(p_2)-h_n^{(p_2)})}{(n-1/2)^q}-\t(p_2+q,p_1)\\
&=-\t(p_2,q,p_1)+\t(p_2)\t(q,p_1)-\t(p_2+q,p_1).
\end{align*}
Hence, from Theorems \ref{thm2.5}-\ref{thm2.8}, we obtain the conclusion that the Hoffman $t$-values of depth three can be expressed in terms of zeta values, double zeta values and double $t$-values.

\subsection{Cubic and Higher Order Euler $T$-sums}

For higher degree sums, like the cubic
\[T_{1^3,q}:=\sum\limits_{n=1}^\infty \frac{h_{n-1}^3}{(n-1/2)^q}\]
it is natural to consider the kernel $\pi \tan(\pi s)\left( \Psi(1/2-s)\right)^3$.

\begin{thm}\label{thm6} For positive integer $q>1$,
\begin{align}\label{2.16}
(1+(-1)^q)T_{1^3,q}&=3\pi^2(\z(q,1)+2\log(2)\z(q))-(q-3)\pi^2\z(q+1)-3T_{1^2,q+1}\nonumber\\ &\quad-3\t(q+2,1)-\t(q+3)+{\rm Res}[F_{0^3,q}(s),s=0],
\end{align}
where $F_{0^3,q}(s)$ defined by
\[F_{0^3,q}(s):=\frac{\pi \tan(\pi s)\left( \Psi(1/2-s)\right)^3}{s^q},\]
and
\begin{align}\label{2.17}
{\rm Res}[F_{0^3,q}(s),s=0]&=(1-(-1)^q)\t(q+3)+(1+(-1)^q)6\log(2)\t(q+2)\nonumber\\&\quad+(1-(-1)^q)12\log^2(2)\t(q+1)+(1+(-1)^q)8\log^3(2)\t(q)\nonumber\\
&\quad-6\sum\limits_{2k_1+k_2=q+2,\atop k_1,k_2\geq 1}\t(2k_1)\z(k_2+1)-24\log(2)\sum\limits_{2k_1+k_2=q+1,\atop k_1,k_2\geq 1}\t(2k_1)\z(k_2+1)\nonumber\\
&\quad-24\log^2(2)\sum\limits_{2k_1+k_2=q,\atop k_1,k_2\geq 1}\t(2k_1)\z(k_2+1)\nonumber\\&\quad+6\sum\limits_{2k_1+k_2+k_3=q+1,\atop k_1,k_2,k_3\geq 1}\t(2k_1)\z(k_2+1)\z(k_3+1)\nonumber\\
&\quad+12\log(2)\sum\limits_{2k_1+k_2+k_3=q,\atop k_1,k_2,k_3\geq 1}\t(2k_1)\z(k_2+1)\z(k_3+1)\nonumber\\
&\quad-2\sum\limits_{2k_1+k_2+k_3+k_4=q,\atop k_1,k_2,k_3\geq 1}\t(2k_1)\z(k_2+1)\z(k_3+1)\z(k_4+1).
\end{align}
\end{thm}
\pf The function $F_{0^3,q}(s)$ only have poles at the $s=0,\pm (n-1/2),n\ (n\in\N)$. For positive integer $n$, from (\ref{2.2}) we have
\begin{align*}
\left(\Psi(1/2-s)+\gamma \right)^3=\frac{1}{(s-n)^3}+3\frac{H_n+2\log(2)}{(s-n)^2}+3\frac{(H_n+2\log(2))^2-(H_n^{(2)}+\z(2))}{s-n}+\cdots,
\end{align*}
Hence, by Lemma \ref{lem2} and (\ref{2.11}),
\begin{align*}
{\rm Res}[F_{0^3,q}(s),s=n]&=\lim\limits_{s\rightarrow n}\frac{1}{2!} \frac{d^2}{ds^2} \left[\left(1+3(H_n+2\log(2))(s-n)+\cdots\right)\frac{\pi \tan(\pi s)}{s^q}\right]\\
&=-q\frac{\pi^2}{n^{q+1}}+3\pi^2\frac{H_n+2\log(2)}{n^q}\quad(n\in\N).
\end{align*}
Moreover,
\begin{align*}
&{\rm Res}[F_{0^3,q}(s),s=n-1/2]=-\frac{h_n^3}{(n-1/2)^q}\quad (n\in\N),\\
&{\rm Res}[F_{0^3,q}(s),s=1/2-n]=(-1)^{q+1}\frac{h_{n-1}^3}{(n-1/2)^q}\quad(n\in\N).
\end{align*}
If $s\rightarrow 0$, then
\begin{align*}
F_{0^3,q}(s)=\frac{2}{s^q}\sum\limits_{k=1}^\infty \t(2k)s^{2k-1} \left( \frac{1}{s}+2\log(2)-\sum\limits_{j=1}^\infty \z(j+1)s^j\right)^3.
\end{align*}
By direct calculations, we obtain the formula (\ref{2.17}). Hence, Theorem \ref{thm6} holds.\hfill$\square$

As an example, setting $q=2$ in (\ref{2.16}) yields
\[T_{1^3,2}+\frac{3}{2}T_{1^2,3}=\frac{5}{4}\pi^2\z(3)+4\log^3(2)\pi^2+\frac{31}{4}\z(5).\]
In general, consider the function
\begin{align*}
F_{(p_1-1)(p_2-1)\cdots(p_k-1),q}(s):=\frac{\pi\tan(\pi s)\Psi^{(p_1-1)}(1/2-s)\Psi^{(p_2-1)}(1/2-s)\cdots \Psi^{(p_k-1)}(1/2-s)}{s^q(p_1-1)!(p_2-1)!\cdots(p_k-1)!}
\end{align*}
with an elementary calculation, we can get the following theorem.
\begin{thm}
A nonlinear Euler T-sum $T_{p_1p_2\cdots p_k,q}$ reduces to a combination of $\log(2)$, Euler T-sums with degree $\leq k-1$, multiple zeta values with depth $\leq k$ whenever the weight $p_1+p_2+\cdots+p_k+q$ and the order $k$ are of the same parity.
\end{thm}

\section{Evaluations of Alternating Euler $T$-sums}

We define alternating Euler $T$-sums ${\bar T}_{p_1\cdots p_k,q}$ by
\begin{align}
{\bar T}_{p_1p_2\cdots p_k,q}:=\sum\limits_{n=1}^\infty \frac{h_{n-1}^{(p_1)}h_{n-1}^{(p_2)}\cdots h_{n-1}^{(p_k)}}{(n-1/2)^q}(-1)^n.
\end{align}
In (\ref{1.1})-(\ref{1.3}), we put a bar on top of $s_j\ (j=1,\cdots k)$ if there is a sign $(-1)^{n_j}$ appearing in the denominator on the right, which are called the alternating MZVs, alternating multiple $t$-values and multiple $\t$-values, respectively. For example,
\[\z({\bar s_1},s_2,{\bar s_3})=\sum\limits_{n_1>n_2>n_3\geq 1}\frac{(-1)^{n_1+n_3}}{n_1^{s_1}n_2^{s_2}n_3^{s_3}},\quad \t({\bar s_1},s_2)=\sum\limits_{n_1>n_2\geq 1}\frac{(-1)^{n_1}}{(n_1-1/2)^{s_1}(n_2-1/2)^{s_2}}.\]
Some results of alternating MZVs, readers may consult \cite{BBV2010,BBBL1997,BBBL2001,BB2003,M2000} and references therein.

For convenience, we let
\[{\bar t}(s):=-\t({\bar s})=\sum\limits_{n=1}^\infty \frac{(-1)^{n-1}}{(n-1/2)^s}.\]
In this section, we will discuss the alternating Euler $T$-sums and find some evaluations of closed form of it. We need the formula (see \cite{XC2019})
\begin{align}\label{3.2}
\frac{\pi }
{{\cos \left( {\pi s} \right)}}\mathop  = \limits^{s \to n-1/2} {\left( { - 1} \right)^n}\left\{ {\frac{1}
{{s - \frac{{2n - 1}}
{2}}} - 2\sum\limits_{k = 1}^\infty  { \zeta \left( {\overline {2k}} \right){{\left( {s - \frac{{2n - 1}}
{2}} \right)}^{2k - 1}}} } \right\}.
\end{align}
From \cite{G1989},
\begin{align}\label{3.3}
\frac{\pi}{\cos(\pi s)}=2\sum\limits_{k=0}^{\infty} {\bar t}(2k+1)s^{2k}=\sum\limits_{k=0}^{\infty} \frac{(-1)^kE_{2k}\pi^{2k+1}}{(2k)!}s^{2k},
\end{align}
where $E_{2k}$ is Euler number. Then, we compute
\begin{align}\label{3.4}
\lim\limits_{s\rightarrow n} \frac{d^p}{ds^p}\frac{\pi }{\cos(\pi s)}=(-1)^n(1+(-1)^p)p!{\bar t}(p+1).
\end{align}

\subsection{Linear Alternating Euler $T$-sums}

\begin{thm}\label{thm3.1} For positive integer $q$,
\begin{align}\label{3.5}
(1+(-1)^{q+1}){\bar T}_{1,q}&=(-1)^{q+1}{\bar t}(q+1)-\pi \z({\bar q})+((-1)^q-1)2\log(2){\bar t}(q)\nonumber\\&\quad+2\sum\limits_{2k_1+k_2=q-2,\atop k_1,k_2\geq1} {\bar t}(2k_1+1)\z(k_2+2).
\end{align}
\end{thm}

\begin{thm}\label{thm3.2} For positive integer $p>1$ and $q$,
\begin{align}\label{3.6}
(1+(-1)^{p+q}){\bar T}_{p,q}&=(-1)^{p+q}{\bar t}(p+q)+(-1)^{p}(1-(-1)^q)\t(p){\bar t}(q)\nonumber\\
&\quad+(-1)^p\sum\limits_{k=0}^{p-1} ((-1)^k+1)\binom{p+q-k-2}{q-1}{\bar t}(k+1)\z(\overline{p+q-k-1})\nonumber\\
&\quad-2(-1)^p\sum\limits_{2k_1+k_2=q-1,\atop k_1,k_2\geq 0}\binom{p+k_2-1}{p-1}{\bar t}(2k_1+1)\z(p+k_2).
\end{align}
\end{thm}
\pf The proofs of Theorems \ref{thm3.1}-\ref{thm3.2} are similar as the proofs of Theorems \ref{thm1}-\ref{thm2}. We only prove the Theorem \ref{thm3.2}. Consider the function
\begin{align*}
G_{p-1,q}(s):=\frac{\pi\Psi^{(p-1)}(1/2-s)}{\cos(\pi s)s^q(p-1)!}.
\end{align*}
By (\ref{2.5})-(\ref{2.7}), (\ref{3.2})-(\ref{3.4}), these residues are ($n\in\N$)
\begin{align*}
&{\rm Res}[G_{p-1,q}(s),s=n]=(-1)^{p-1}\sum\limits_{k=0}^{p-1} ((-1)^k+1)\binom{p+q-k-2}{q-1}{\bar t}(k+1)\frac{(-1)^n}{n^{p+q-k-1}},\\
&{\rm Res}[G_{p-1,q}(s),s=n-1/2]=\frac{(-1)^p\t(p)+h_n^{(p)}}{(n-1/2)^q}(-1)^n,\\
&{\rm Res}[G_{p-1,q}(s),s=1/2-n]=-\frac{\t(p)-h_{n-1}^{(p)}}{(n-1/2)^q}(-1)^{n+p+q},\\
&\begin{aligned}{\rm Res}[G_{p-1,q}(s),s=0]&=(1-(-1)^{p+q}){\bar t}(p+q)\\&\quad+2(-1)^p\sum\limits_{2k_1+k_2=q-1,\atop k_1,k_2\geq 0}\binom{p+k_2-1}{p-1}{\bar t}(2k_1+1)\z(p+k_2).\end{aligned}
\end{align*}
Thus, by Lemma \ref{lem1} and summing these four contributions yields the desired result.\hfill$\square$

Letting $q=1$ in (\ref{3.5}) and $p=q=2$ in (\ref{3.6}) give
\begin{align*}
&{\bar T}_{1,1}=\t({\bar 1},1)=2G-\frac{1}{2}\pi \log(2),\\
&{\bar T}_{2,2}=\t({\bar 2},2)=\frac{1}{2} \t(4)-\frac{7}{4}\pi \z(3),
\end{align*}
where ${\bar t}(2)=4G$, $G$ is Catalan's constant.

\subsection{Quadratic Alternating Euler $T$-sums}

\begin{thm} For positive integer $q>1$,
\begin{align}\label{3.7}
(1+(-1)^{q+1}){\bar T}_{1^2,q}&=(q-2)\pi \z(\overline{q+1})-2\pi \z({\bar q},1)-4\pi \log(2)\z(\bar q)-2\t(\overline{q+1},1)\nonumber\\
&\quad+(-1)^q{\bar t}(q+2)-(1+(-1)^q)4\log(2){\bar t}(q+1)-(1-(-1)^q) 4\log^2(2){\bar t}(q)\nonumber\\
&\quad-2\sum\limits_{2k_1+k_2+k_3=q+1,\atop k_1,k_2,k_3\geq 1} {\bar t}(2k_1-1)\z(k_2+1)\z(k_3+1)\nonumber\\&\quad+4\sum\limits_{2k_1+k_2=q+2,\atop k_2,k_2\geq 1}{\bar t}(2k_1-1)\z(k_2+1)
\nonumber\\&\quad+8\log(2)\sum\limits_{2k_1+k_2=q+1,\atop k_2,k_2\geq 1}{\bar t}(2k_1-1)\z(k_2+1).
\end{align}
\end{thm}
\pf The proof is based on the function
\begin{align*}
G_{0^2,q}(s):=\frac{\pi\left(\Psi(1/2-s)\right)^2}{\cos(\pi s)s^q}
\end{align*}
and the usual residue computation. By a similar argument as in the proof of the above theorem, we may easily deduce the (\ref{3.7}).

If letting $q=1,3$ in (\ref{3.7}), then
\begin{align*}
&{\bar T}_{1^2,1}=\frac{\pi^3}{12}-\frac{\pi }{2}\log^2(2)-\t({\bar 2},1),\\
&{\bar T}_{1^2,3}=\frac{17\pi^5}{720}+\frac{7}{2}\pi \log(2)\z(3)-\log^2(2)\pi^3-\pi\z({\bar 3},1)-\t({\bar 4},1).
\end{align*}

Similarly, by considering the function
\begin{align*}
G_{(p_1-1)(p_2-1)\cdots(p_k-1),q}(s):=\frac{\pi\Psi^{(p_1-1)}(1/2-s)\Psi^{(p_2-1)}(1/2-s)\cdots \Psi^{(p_k-1)}(1/2-s)}{\cos(\pi s)s^q(p_1-1)!(p_2-1)!\cdots(p_k-1)!}
\end{align*}
and using the residue computations, many other relations involving alternating Euler $T$-sums can be established.

\section{Other Euler Type $T$-Sums}
In this section, we define an Euler type $T$-sums
\begin{align*}
\S_{p_1p_2\cdots p_k,q}:=\sum\limits_{n=1}^\infty \frac{h_n^{(p_1)}h_n^{(p_2)}\cdots h_n^{(p_k)}}{n^q}.
\end{align*}
Then, we apply the contour integral to establish many relations of $\S_{p_1p_2\cdots p_k,q}$. Moreover, we can use the quadratic sum $\S_{p_1p_2,q}$ to evaluate the triple Kaneko-Tsumura $T$-values with weight even.

The Kaneko-Tsumura multiple $T$-zeta values are defined by \cite{KTA2018,KTA2019}
\begin{align}
T(k_1,k_2,\ldots,k_r):=2^r\sum\limits_{n_1>n_2>\cdots>n_r>0} \frac{1}{(2n_1-r)^{k_1}(2n_2-r+1)^{k_2}\cdots (2n_r-1)^{k_r}},
\end{align}
where we used the opposite convention ($n_1>n_2>\cdots>n_r>0$) of the original definition ($0<n_1<\cdots<n_r$) Kaneko-Tsumura $T$-zeta values ($k_j\in \N,\ k_1>1$). Here $k_1+k_2+\cdots+k_r$ and $r$ are called the weight and depth of the sum, respectively.

Hence, from the definitions of $S_{p_1p_2\cdots p_k,q}$ and $T$-values, we have
\begin{align}
&T(k_1,k_2)=\frac{1}{2^{k_1+k_2-2}}\S_{k_2,k_1},\label{4.3}\\
&T(k_1,k_2,k_3)=\frac{1}{2^{k_1+k_2+k_3-3}}\t(k_1)\S_{k_3,k_2}-\frac{1}{2^{k_1+k_2+k_3-3}}\S_{k_1k_3,k_2}.\label{4.4}
\end{align}
On the other hand, according to the relation of double $T$-values and double zeta values,
\begin{align*}
T(k_1,k_2)=\zeta({\bar k}_1,k_2)+\zeta(k_1,k_2)-\zeta({k_1},{\bar k}_2)-\zeta({\bar k_1},{\bar k_2})
\end{align*}
we obtain
\begin{align}\label{4.5}
T(q,p)=&\frac{1-(-1)^p}{2} (\zeta(p)-\zeta(\bar p))\nonumber\\&+(-1)^p \sum\limits_{k=0}^{[p/2]} \binom{m-2k-1}{q-1} (\zeta(2k)-\zeta(\overline{2k}))(\zeta(m-2k)-\zeta(\overline{m-2k}))\nonumber\\
&+(-1)^p \sum\limits_{k=0}^{[q/2]} \binom{m-2k-1}{p-1} (\zeta(2k)+\zeta(\overline{2k}))(\zeta(m-2k)-\zeta(\overline{m-2k})),
\end{align}
where $\zeta \left(1\right)$ should be interpreted as $0$ wherever it occurs, and $\zeta \left(0\right)={\zeta}\left(\bar 0\right)=-1/2$, $m=p+q$ odd.

From (\ref{4.3}) and (\ref{4.5}), we know that the linear sum $\S_{p,q}$ with $p+q$ odd can be evaluated by zeta values. In fact, considering the function
\begin{align*}
F_{p-1,q}(s):=\frac{ \Psi^{(p-1)}(1/2-s)}{(s+1/2)^q(p-1)!}\pi\tan(\pi s)
\end{align*}
and using the residue theorem, we can also obtain the result (\ref{4.5}).

By harmonic product, we give $(k_1,k_2\geq 2)$
\begin{align}\label{d5}
\sum\limits_{n=1}^\infty \frac{H_{n-1}^{(k_1)}}{(n-1/2)^{k_2}}=-2^{k_1+k_2-2}T(k_1,k_2)+\z(k_1)\t(k_2).
\end{align}
In fact, from \cite{XC2019}, we also obtain ($p>1$)
\begin{align}
\sum\limits_{n=1}^\infty \frac{H_{n-1}}{(n-1/2)^{p}}=\frac{p}{2}\t(p+1)-\frac{1}{2}\sum\limits_{j=1}^{p-2}\t(p-j)\t(j+1)-2\log(2)\t(p).
\end{align}

For nonnegative integers $n_1,\ldots,n_p$ and $n$, let
\[\binom{n}{n_1,\ldots,n_{p}}:=\frac{n!}{n_1!n_2!\cdots n_p!},\quad (0\leq n_1+\cdots+n_p\leq n).\]
\begin{thm}\label{thm4.1} For positive integers $q>1$,
\begin{align}
(1+(-1)^q)\S_{1^p,q}&=(-1)^{p-1}\sum\limits_{l=0}^{p-1}((-1)^l-1)\binom{p+q-l-2}{q-1}\t(l+1)\t(p+q-l-1)\nonumber\\
&\quad+\sum\limits_{k=1}^p \binom{p}{k} \sum\limits_{|k|_{p-1}=k,\atop 0\leq |{\widetilde{k}}|_{p-1}\leq p-1.} \binom{k}{k_1,\ldots,k_{p-1}} (-1)^{p-1-|{\widetilde{k}}|_{p-1}}\sum\limits_{l=0}^{p-1-|{\widetilde{k}}|_{p-1}} ((-1)^l-1)\nonumber\\
&\quad\quad\quad\quad\quad\quad\times\binom{p+q-|{\widetilde{k}}|_{p-1}-l-2}{q-1}\t(l+1)\sum\limits_{n=1}^\infty \frac{\prod\limits_{j=1}^{p-1}C_{n-1}^{k_j}(j)}{(n-1/2)^{p+q-|{\widetilde{k}}|_{p-1}-l-1}}\nonumber\\
&\quad+(-1)^{p+1}\sum\limits_{k_1+k_2+\cdots+k_p=q,\atop k_j\in\N}\t(k_1+1)\t(k_2+1)\cdots\t(k_p+1)\nonumber\\
&\quad-2(-1)^{p+1}\sum\limits_{2k_1+k_2+\cdots+k_{p+1}=q,\atop k_j\in\N}\z(2k_1)\t(k_2+1)\cdots\t(k_{p+1}+1),
\end{align}
where $|k|_{p-1}:=k_1+k_2+\cdots+k_{p-1}\ (k_j\in\N_0),\ |{\widetilde{k}}|_{p-1}=k_1+2k_2+\cdots+(p-1)k_{p-1}\ (k_j\in\N_0)$ and $\t(1):=0$,
\begin{align}
C_n(j)= \left\{ {\begin{array}{*{20}{c}}\ \ H_n+2\log(2)
   {,\ \ \quad j=1}  \\
   (-1)^{j-1}H_n^{(j)}-\z(j) ,\ j\geq 2. \\
\end{array} } \right.
\end{align}
\end{thm}
\pf The proof is based on the function
\[{\bar F}_{0^p,q}(s):=\frac{\pi \tan(\pi s)\left( \Psi(1/2-s)\right)^p}{(s+1/2)^q},\]
and the usual residue computation. we achieve the desired expansion after a rather tedious computation. \hfill$\square$

If $p=2,3$, then we can get the following corollaries
\begin{cor}  For positive integers $q>1$,
\begin{align}\label{d1}
(1+(-1)^q)\S_{1^2,q}&=\pi^2\t(q)-\sum\limits_{k_1+k_2=q,\atop k_1,k_2\geq 1} \t(k_1+1)\t(k_2+1)\nonumber\\&\quad +2\sum\limits_{2k_1+k_2+k_3=q,\atop k_1,k_2\geq 1} \z(2k_1)\t(k_2+1)\t(k_3+1).
\end{align}
\end{cor}

\begin{cor}  For positive integers $q>1$,
\begin{align}\label{d2}
(1+(-1)^q)\S_{1^3,q}&=-2q\t(2)\t(q+1)+6\t(2)\sum\limits_{n=1}^\infty \frac{H_{n-1}}{(n-1/2)^q}+12\log(2)\t(2)\t(q)\nonumber\\
&\quad+\sum\limits_{k_1+k_2+k_3=q,\atop k_1,k_2,k_3\geq1} \t(k_1+1)\t(k_2+1)\t(k_3+1)
\nonumber\\
&\quad-2\sum\limits_{2k_1+k_2+k_3+k_4=q,\atop k_1,k_2,k_3,k_4\geq1} \z(2k_1)\t(k_2+1)\t(k_3+1)\t(k_4+1).
\end{align}
\end{cor}

Setting $q=2$ and $4$ in (\ref{d1}) and (\ref{d2}) yield
\begin{align*}
&\S_{1^2,2}=\frac{\pi^4}{8},\quad \S_{1^2,4}=\frac{\pi^6}{24}-\frac{49}{2}\zeta^2(3),\\
&\S_{1^3,2}=\frac7{2}\pi^2\z(3),\quad \S_{1^3,4}=-\frac{21}{8}\pi^4\z(3)+31\pi^2\z(5).
\end{align*}

\begin{thm}\label{thm4.4} For positive integers $p,q>1$,
\begin{align}
(1-(-1)^{p+q})\S_{1p,q}&=(-1)^{p-1}(1+(-1)^q)\t(p)\S_{1,q}\nonumber\\&\quad+(-1)^p\sum\limits_{l=0}^p ((-1)^l-1)\binom{p+q-l-1}{q-1}\t(l+1)\t(p+q-l)\nonumber\\&\quad
+(-1)^{p-1}\sum\limits_{l=0}^{p-1}((-1)^l-1)\binom{p+q-l-2}{q-1}\t(l+1)\nonumber\\&\quad\quad\quad\quad\times\left(\sum\limits_{n=1}^\infty \frac{H_{n-1}}{(n-1/2)^{p+q-l-1}}+2\log(2)\t(p+q-l-1) \right)\nonumber\\&\quad
+(-1)^{p-1}\sum\limits_{k=1}^{p-1}\sum\limits_{l=0}^{p-k-1}((-1)^l-1)\binom{p+q-k-l-2}{q-1}\t(l+1)\nonumber\\&\quad
\quad\quad\quad\times\left(\sum\limits_{n=1}^\infty \frac{H_{n-1}^{(k+1)}}{(n-1/2)^{p+q-k-l-1}} -(-1)^k \z(k+1)\t(p+q-k-l-1)\right)\nonumber\\&\quad
+(-1)^p\sum\limits_{k_1+k_2=q+1,\atop k_1,k_2\geq 1} \binom{k_2+p-2}{p-1}\t(k_1+1)\t(k_2+p-1)\nonumber\\&\quad
-2(-1)^p\sum\limits_{2k_1+k_2+k_3=q+1,\atop k_1,k_2,k_3\geq 1} \binom{k_3+p-2}{p-1}\z(2k_1)\t(k_2+1)\t(k_3+p-1).
\end{align}
\end{thm}
\pf By considering the function
\[{\bar F}_{0(p-1),q}(s):=\frac{\pi \tan(\pi s) \Psi(1/2-s)\Psi^{(p-1)}(1/2-s)}{(s+1/2)^q(p-1)!},\]
and using the direct residue computation. The desired formulas can be established. \hfill$\square$

If $p=2,q=3$, then
\[\S_{12,3}=-\frac{\pi^6}{16}+49\zeta^2(3).\]

\begin{thm}\label{thm4.5} For positive integers $p_1,p_2,q>1$,
\begin{align}
&(1+(-1)^{p_1+p_2+q})\S_{p_1p_2,q}\nonumber\\&=-(-1)^{p_1+p_2}(1+(-1)^q)\t(p_1)\t(p_2)\z(q)\nonumber\\&\quad+(-1)^{p_1}((-1)^{p_2+q}-1)\t(p_1)\S_{p_2,q}+(-1)^{p_2}((-1)^{p_1+q}-1)\t(p_2)\S_{p_1,q}
\nonumber\\&\quad-(-1)^{p_1+p_2}\sum\limits_{k_1+k_2=q+2,\atop k_1,k_2\geq 1} \binom{k_1+p_1-2}{p_1-1}\binom{k_2+p_2-2}{p_2-1}\t(k_1+p_1-1)\t(k_2+p_2-1)\nonumber\\&\quad
+2(-1)^{p_1+p_2}\sum\limits_{k_1+k_2+2k_3=q+2,\atop k_1,k_2,k_3\geq 1}\binom{k_1+p_1-2}{p_1-1}\binom{k_2+p_2-2}{p_2-1}\t(k_1+p_1-1)\t(k_2+p_2-1)\z(2k_3)\nonumber\\&\quad
-(-1)^{p_1+p_2}\sum\limits_{l=0}^{p_1+p_2-1}((-1)^l-1)\binom{p_1+p_2+q-l-2}{q-1}\t(l+1)\t(p_1+p_2+q-l-1)\nonumber\\&\quad
+(-1)^{p_1+p_2}\sum\limits_{k=1}^{p_2}(-1)^k\binom{k+p_1-2}{p_1-1}\sum\limits_{l=0}^{p_2-k}((-1)^l-1)\binom{p_2+q-k-l-1}{q-1}\t(l+1)\nonumber\\&\quad
\quad\quad\times\left(\z(k+p_1-1)\t(p_2+q-k-l)+(-1)^{k+p_1-1}\sum\limits_{n=1}^\infty \frac{H_{n-1}^{(k+p_1-1)}}{(n-1/2)^{p_2+q-k-l}} \right)\nonumber\\&\quad
+(-1)^{p_1+p_2}\sum\limits_{k=1}^{p_1} (-1)^k\binom{k+p_2-2}{p_2-1}\sum\limits_{l=0}^{p_1-k}((-1)^l-1)\binom{p_1+q-k-l-1}{q-1}\t(l+1)\nonumber\\&\quad
\quad\quad\times\left(\z(k+p_2-1)\t(p_1+q-k-l)+(-1)^{k+p_2-1}\sum\limits_{n=1}^\infty \frac{H_{n-1}^{(k+p_2-1)}}{(n-1/2)^{p_1+q-k-l}} \right).
\end{align}
\end{thm}
\pf By considering the function
\[{\bar F}_{(p_1-1)(p_2-1),q}(s):=\frac{\pi \tan(\pi s)\Psi^{(p_1-1)}(1/2-s)\Psi^{(p_2-1)}(1/2-s)}{(s+1/2)^q(p_1-1)!(p_2-1)!},\]
then a direct residue computation gives the desired formula.\hfill$\square$

Hence, from the relations (\ref{4.3}), (\ref{4.4}) and Theorems \ref{thm4.1}, \ref{thm4.4}, \ref{thm4.5}, we obtain the formula of the triple $T$-value $T(p, q, r)$ of even weight in terms of single and the double $T$-values. Tsumura also proved an explicit formula of triple $T$-value, see \cite{T2019} for the detail.

Setting $p_1=p_2=q=2$ yields
\begin{align*}
\S_{2^2,2}&=\t(2)\t(4)-2\t^2(3)+2\t^2(2)\z(2)+2\t(2)\sum\limits_{n=1}^\infty \frac{H_{n-1}^{(2)}}{(n-1/2)^2}\\&
=32 \pi ^2 \text{Li}_4\left(\frac{1}{2}\right)-98 \zeta^2 (3)+28 \pi ^2 \zeta (3) \log (2)-\frac{61 \pi ^6}{360}+\frac{4}{3} \pi ^2 \log ^4(2)-\frac{4}{3} \pi ^4 \log ^2(2).
\end{align*}

In general, we can consider the function
\begin{align*}
{\bar F}_{(p_1-1)(p_2-1)\cdots(p_k-1),q}(s):=\frac{\pi\tan(\pi s)\Psi^{(p_1-1)}(1/2-s)\Psi^{(p_2-1)}(1/2-s)\cdots \Psi^{(p_k-1)}(1/2-s)}{(s+1/2)^q(p_1-1)!(p_2-1)!\cdots(p_k-1)!}
\end{align*}
to establish more general formulas of $\S_{p_1p_2\cdots p_k,q}$, but it is very difficult. Similarly, we can also consider the function
\begin{align*}
{\bar G}_{(p_1-1)(p_2-1)\cdots(p_k-1),q}(s):=\frac{\pi\Psi^{(p_1-1)}(1/2-s)\Psi^{(p_2-1)}(1/2-s)\cdots \Psi^{(p_k-1)}(1/2-s)}{\cos(\pi s)(s+1/2)^q(p_1-1)!(p_2-1)!\cdots(p_k-1)!}
\end{align*}
to evaluate the alternating sum $\S_{p_1p_2\cdots p_k,{\bar q}}$ defined by

\begin{align*}
\S_{p_1p_2\cdots p_k,{\bar q}}:=\sum\limits_{n=1}^\infty \frac{h_n^{(p_1)}h_n^{(p_2)}\cdots h_n^{(p_k)}}{n^q}(-1)^n.
\end{align*}
It is possible that of some other relations involving alternating Euler $T$-sums can be proved by using the
techniques of the present paper. For example, we can define an alternating ${\bar \Psi}(-s)$ function
$${\bar \Psi}(-s):=\frac{1}{s+1/2}+\sum\limits_{k=1}^\infty \left(\frac{(-1)^k}{k-1/2}-\frac{(-1)^k}{k-1/2-s} \right).$$
Then, consider the four functions
\begin{align*}
&E_{(p_1-1)(p_2-1)\cdots(p_k-1),q}(s):=\frac{\pi\tan(\pi s){\bar \Psi}^{(p_1-1)}(1/2-s){\bar \Psi}^{(p_2-1)}(1/2-s)\cdots {\bar \Psi}^{(p_k-1)}(1/2-s)}{s^q(p_1-1)!(p_2-1)!\cdots(p_k-1)!},\\
&H_{(p_1-1)(p_2-1)\cdots(p_k-1),q}(s):=\frac{\pi{\bar \Psi}^{(p_1-1)}(1/2-s){\bar \Psi}^{(p_2-1)}(1/2-s)\cdots {\bar \Psi}^{(p_k-1)}(1/2-s)}{\cos(\pi s)s^q(p_1-1)!(p_2-1)!\cdots(p_k-1)!},\\
&{\bar E}_{(p_1-1)(p_2-1)\cdots(p_k-1),q}(s):=\frac{\pi\tan(\pi s){\bar \Psi}^{(p_1-1)}(1/2-s){\bar \Psi}^{(p_2-1)}(1/2-s)\cdots {\bar \Psi}^{(p_k-1)}(1/2-s)}{(s+1/2)^q(p_1-1)!(p_2-1)!\cdots(p_k-1)!},\\
&{\bar H}_{(p_1-1)(p_2-1)\cdots(p_k-1),q}(s):=\frac{\pi {\bar \Psi}^{(p_1-1)}(1/2-s){\bar \Psi}^{(p_2-1)}(1/2-s)\cdots {\bar \Psi}^{(p_k-1)}(1/2-s)}{\cos(\pi s)(s+1/2)^q(p_1-1)!(p_2-1)!\cdots(p_k-1)!},
\end{align*}
and use the contour integral to evaluate these sums

\begin{align*}
&\sum\limits_{n=1}^\infty \frac{{\bar h}_{n-1}^{(p_1)}{\bar h}_{n-1}^{(p_2)}\cdots {\bar h}_{n-1}^{(p_k)}}{(n-1/2)^q},\quad \sum\limits_{n=1}^\infty \frac{{\bar h}_{n-1}^{(p_1)}{\bar h}_{n-1}^{(p_2)}\cdots {\bar h}_{n-1}^{(p_k)}}{(n-1/2)^q}(-1)^n,\\
&\sum\limits_{n=1}^\infty \frac{{\bar h}_n^{(p_1)}{\bar h}_n^{(p_2)}\cdots {\bar h}_n^{(p_k)}}{n^q},\quad \sum\limits_{n=1}^\infty \frac{{\bar h}_n^{(p_1)}{\bar h}_n^{(p_2)}\cdots {\bar h}_n^{(p_k)}}{n^q}(-1)^n
\end{align*}
where ${\bar h}_n^{(p)}$ is defined by
\begin{align*}
{\bar h}_n^{(p)}:=\sum\limits_{k=1}^n \frac{(-1)^n}{(n-1/2)^p},\quad {\bar h}_0^{(p)}:=0.
\end{align*}

\section{Formulas of Kaneko-Tsumura's Conjecture}

In \cite{KTA2019}, Kaneko and Tsumura conjecture the following relation $(p\geq 2,m,q\geq 1)$
\begin{align*}
\sum\limits_{i+j=m,\atop i,j\geq 0} \binom{p+i-1}{i}\binom{q+j-1}{j} T(p+i,q+j)\in \ZZ
\end{align*}
where $\ZZ$ is the space of usual multiple zeta values. Quite recently, T. Murakami proved the conjecture by using the motivic method employed in \cite{Gla}, but not gave explicit formula. In this section, we will give an explicit duality formula of the conjecture by using the residue theorem.

\begin{thm}\label{thm5.1} For positive integers $m,p$ and $q>1$,
\begin{align}\label{5.1}
&(-1)^{m-1} \sum\limits_{i+j=p-1,\atop i,j\geq 0} \binom{m+i-1}{i}\binom{q+j-1}{j} \sum\limits_{n=1}^\infty \frac{H_{n-1}^{(m+i)}}{(n-1/2)^{q+j}}\nonumber\\
+&(-1)^{p-1} \sum\limits_{i+j=m-1,\atop i,j\geq 0} \binom{p+i-1}{i}\binom{q+j-1}{j} \sum\limits_{n=1}^\infty \frac{H_{n-1}^{(p+i)}}{(n-1/2)^{q+j}}\nonumber\\
&=\binom{p+q+m-2}{q-1}\t(p+q+m-1) \nonumber\\&\quad+ \sum\limits_{i+j=p-1,\atop i,j\geq 0} \binom{m+i-1}{i}\binom{q+j-1}{j} (-1)^i \z(m+i)\t(q+j)\nonumber\\&\quad+\sum\limits_{i+j=m-1,\atop i,j\geq 0} \binom{p+i-1}{i}\binom{q+j-1}{j} (-1)^i \z(p+i)\t(q+j)\nonumber\\
&\quad-\sum\limits_{i+j=q-1,\atop i,j\geq 0} \binom{m+i-1}{i}\binom{p+j-1}{j} \t(m+i)\t(p+j),
\end{align}
where $\zeta \left(1\right)$ should be interpreted as $-2\log(2)$ wherever it occurs, and $\t(1):=0$.
\end{thm}
\pf To prove the identity, we consider the function
\begin{align*}
F(s):=\frac{\Psi^{(m-1)}(1/2-s)\Psi^{(p-1)}(1/2-s)}{(s+1/2)^q(m-1)!(p-1)!}.
\end{align*}
It is obvious that the function $F(s)$ only have poles at the $s=-1/2$ and $s=n\ (n\in\N_0)$. At a non-negative integer $n$, by (\ref{2.2}) and (\ref{2.5}), we compute the residue
\begin{align*}
{\rm Res}[F(s),s=n]=&-(-1)^{m+p}\binom{p+q+m-2}{q-1} \frac{1}{(n+1/2)^{p+q+m-1}}\\
&-(-1)^{m+p} \sum\limits_{i+j=p-1,\atop i,j\geq 0} \binom{m+i-1}{i}\binom{q+j-1}{j} \frac{(-1)^{i}\z(m+i)+(-1)^{m}H_n^{(m+i)}}{(n+1/2)^{q+j}}\\
&-(-1)^{m+p} \sum\limits_{i+j=m-1,\atop i,j\geq 0} \binom{p+i-1}{i}\binom{q+j-1}{j} \frac{(-1)^{i}\z(p+i)+(-1)^{p}H_n^{(m+i)}}{(n+1/2)^{q+j}}.
\end{align*}
By (\ref{2.3}) and (\ref{2.6}), the residue of the pole of order $q$ is found to be
\begin{align*}
{\rm Res}[F(s),s=-1/2]=(-1)^{m+p}\sum\limits_{i+j=q-1,\atop i,j\geq 0} \binom{m+i-1}{i}\binom{p+j-1}{j} \t(m+i)\t(p+j).
\end{align*}
Hence, summing these two contributions yields the desired formula. \hfill$\square$

Letting $p=m=2,q=3$ in (\ref{5.1}), we have
\begin{align*}
3\sum\limits_{n=1}^\infty \frac{H_{n-1}^{(2)}}{(n-1/2)^4}+2\sum\limits_{n=1}^\infty \frac{H_{n-1}^{(3)}}{(n-1/2)^3}=112\zeta^2(3)-\frac{\pi^6}{6}.
\end{align*}

Similarly, by considering the function
\begin{align*}
F(s,a):=\frac{\Psi^{(m-1)}(1/2-s;a)\Psi^{(p-1)}(1/2-s;a)}{(s+1/2)^q(m-1)!(p-1)!}
\end{align*}
with the usual residue computation, we can get the following a more general result.
\begin{thm}\label{thm5.2} For positive integers $m,p$ and $q>1$,
\begin{align}
&(-1)^{m-1} \sum\limits_{i+j=p-1,\atop i,j\geq 0} \binom{m+i-1}{i}\binom{q+j-1}{j} \sum\limits_{n=1}^\infty \frac{H_{n-1}^{(m+i)}}{(n+a)^{q+j}}\nonumber\\
+&(-1)^{p-1} \sum\limits_{i+j=m-1,\atop i,j\geq 0} \binom{p+i-1}{i}\binom{q+j-1}{j} \sum\limits_{n=1}^\infty \frac{H_{n-1}^{(p+i)}}{(n+a)^{q+j}}\nonumber\\
&=\binom{p+q+m-2}{q-1}\z(p+q+m-1;a+1) \nonumber\\&\quad+ \sum\limits_{i+j=p-1,\atop i,j\geq 0} \binom{m+i-1}{i}\binom{q+j-1}{j} (-1)^i \z(m+i)\z(q+j;a+1)\nonumber\\&\quad+\sum\limits_{i+j=m-1,\atop i,j\geq 0} \binom{p+i-1}{i}\binom{q+j-1}{j} (-1)^i \z(p+i)\z(q+j;a+1)\nonumber\\
&\quad-\sum\limits_{i+j=q-1,\atop i,j\geq 0} \binom{m+i-1}{i}\binom{p+j-1}{j} \z(m+i;a+1)\z(p+j;a+1),
\end{align}
where $\zeta \left(1\right)$ should be interpreted as $-2\log(2)$ wherever it occurs, and $\z(1;a+1):=0$. Here $\z(s;a+1)$ stands for Hurwitz zeta function, which is defined by
\[\zeta( {s;a+1}): = \sum\limits_{n = 1}^\infty  {\frac{1}{{{{\left( {n + a} \right)}^s}}}}\quad \left( {{\Re} \left( s \right) > 1,\ a\in \mathbb{C}\setminus \N^-} \right).\]
\end{thm}
It is clear that Theorem \ref{thm5.1} is immediate corollary of Theorem \ref{thm5.2} with $a=-1/2$. Next, we give a duality formulas of Kaneko-Tsumura conjecture.

\begin{thm} For positive integers $p,q,m\geq 2$,
\begin{align}\label{5.2}
&(-1)^{m} \sum\limits_{i+j=p-1,\atop i,j\geq 0} \binom{m+i-1}{i}\binom{q+j-1}{j}T(m+i,q+j)\nonumber\\
+&(-1)^{p} \sum\limits_{i+j=m-1,\atop i,j\geq 0} \binom{p+i-1}{i}\binom{q+j-1}{j} T(p+i,q+j)\nonumber\\
&=\frac{1}{2^{p+q+m-3}}\binom{p+q+m-2}{q-1}\t(p+q+m-1) \nonumber\\&\quad+\frac{1}{2^{p+q+m-3}}\sum\limits_{i+j=p-1,\atop i,j\geq 0} \binom{m+i-1}{i}\binom{q+j-1}{j} \left((-1)^i+(-1)^m \right) \z(m+i)\t(q+j)\nonumber\\&\quad+\frac{1}{2^{p+q+m-3}} \sum\limits_{i+j=m-1,\atop i,j\geq 0} \binom{p+i-1}{i}\binom{q+j-1}{j} \left((-1)^i+(-1)^p \right) \z(p+i)\t(q+j)\nonumber\\
&\quad-\frac{1}{2^{p+q+m-3}}\sum\limits_{i+j=q-1,\atop i,j\geq 0} \binom{m+i-1}{i}\binom{p+j-1}{j} \t(m+i)\t(p+j).
\end{align}
\end{thm}
\pf The result immediately follows from (\ref{d5}) and (\ref{5.1}).   \hfill$\square$

Taking $m=p$ in (\ref{5.2}), we can get the following corollary.
\begin{cor} For positive integers $p,q\geq 2$
\begin{align}\label{5.3}
&\sum\limits_{i+j=p-1,\atop i,j\geq 0} \binom{p+i-1}{i}\binom{q+j-1}{j}T(p+i,q+j)\nonumber\\
&=\frac{1}{2^{2p+q-2}} \binom{2p+q-2}{q-1} \t(2p+q-1)\nonumber\\&\quad+\frac{1}{2^{2p+q-3}}\sum\limits_{i+j=p-1,\atop i,j\geq 0} \binom{p+i-1}{i}\binom{q+j-1}{j} \left(1-(-1)^j\right) \z(p+i)\t(q+j)\nonumber\\
&\quad-\frac{(-1)^p}{2^{2p+q-2}}\sum\limits_{i+j=q-1,\atop i,j\geq 0} \binom{p+i-1}{i}\binom{p+j-1}{j}\t(p+i)\t(p+j).
\end{align}
\end{cor}

For example, setting $p=2,q=3$ in (\ref{5.3}) gives
\[3T(2,4)+2T(3,3)=\frac{\pi^6}{64}-\frac{49}{8}\zeta^2(3).\]

\begin{thm}\label{thm5.4} For positive integers $p,m$ and $q\geq 2$,
\begin{align}\label{5.4}
&(-1)^{p-1} \sum\limits_{i+j=m-1,\atop i,j\geq 0} \binom{p+i-1}{i}\binom{q+j-1}{j}\S_{p+i,q+j}\nonumber\\
&+(-1)^{m-1}\sum\limits_{i+j=p-1,\atop i,j\geq 0} \binom{m+i-1}{i}\binom{q+j-1}{j} T_{m+i,q+j}\nonumber\\
&= \sum\limits_{i+j=m-1,\atop i,j\geq 0} \binom{p+i-1}{i}\binom{q+j-1}{j}(-1)^i \t(p+i)\z(q+j)\nonumber\\
&\quad+ \sum\limits_{i+j=p-1,\atop i,j\geq 0} \binom{m+i-1}{i}\binom{q+j-1}{j}(-1)^i \t(m+i)\t(q+j)\nonumber\\
&\quad- \sum\limits_{i+j=q-1,\atop i,j\geq 0} \binom{m+i-1}{i}\binom{p+j-1}{j}\z(m+i)\t(p+j),
\end{align}
where $\zeta \left(1\right):=-2\log(2)$ and $\t(1):=0$.
\end{thm}
\pf The proof of Theorem \ref{thm5.4} is similar as the proof of Theorem \ref{thm5.1}. We consider the function
\begin{align*}
G(s):=\frac{\Psi^{(m-1)}(1/2-s)\Psi^{(p-1)}(-s)}{(s+1)^q(m-1)!(p-1)!},
\end{align*}
then by a similar argument as in the proof (\ref{5.4}), we deduce the desired result.  \hfill$\square$

Letting $(p,q,m)=(1,2,2)$ and $(2,2,2)$ yield
\begin{align*}
&\S_{2,2}+2\S_{1,3}-T_{2,2}=\frac{\pi^4}{12},\\
&\S_{2,3}+\S_{3,2}+T_{2,3}+T_{3,2}=2\z(2)\t(3).
\end{align*}
Hence, by the relations
\begin{align*}
\S_{p,q}=2^{p+q-2}T(q,p)\quad{\rm and}\quad  T_{p,q}=2^{p+q}t(q,p)
\end{align*}
we can obtain a kind of relationship of the double Kaneko-Tsumura $T$-values and the double Hoffman $t$-values.

{\bf Acknowledgments.}  The authors express their deep gratitude to Professor Masanobu Kaneko and Jianqiang Zhao for valuable discussions and comments. The first author is supported by the China Scholarship Council (No. 201806310063). The second author is supported by the National Natural Science Foundation of China (under Grant 11671360) and the Zhejiang Provincial Natural Science Foundation of China (under Grant LQ17A010010).

 \end{document}